\documentclass{article}
\usepackage{graphicx} 
\usepackage[a4paper, margin=2.5cm]{geometry} 
\usepackage{amsmath}
\usepackage{amssymb}
\usepackage{subcaption}

\title{Exact Autocorrelation of the Multiplicative Inverse of a Non-Zero-Mean Complex Gaussian Process}
\author{Marco Lanucara}

\begin{document}

\maketitle

\begin{abstract}
    We study the spectral properties of a stochastic process obtained by multiplicative inversion of a non-zero-mean complex Gaussian process. 
    We show that its autocorrelation function and power spectrum exist for most regular processes, and we derive a first order differential 
    equation for the autocorrelation function, which follows from a demonstrated recurrence relation between the coefficients of its series 
    expansion. Finally, we obtain a closed-form expression of the autocorrelation function based on the first and second order statistics of 
    the underlying Gaussian process.
\end{abstract}

\section*{Keywords}
Statistical signal processing; stochastic process; Gaussian process.

\section{Introduction}
\label{sec:introduction}

Even though the study of nonlinearities applied to stochastic processes is vast (e.g. \cite{deutsch2017nonlinear}), the case of multiplicative 
inversion seems to have been devoted to zero-mean processes, i.e. to the study of objects like $1/\mathbf{x}(t)$\footnote{Lowercase bold roman 
font is used to identify stochastic processes when an explicit time dependency is reported, in line with \cite{papoulis1965random}, or, 
without time dependency,  4-dimensional vectors and (if uppercase) matrices.}, where $\mathbf{x}(t)$ is a complex zero-mean stochastic 
process. As a related example, the study of the spectral properties of the ratio between complex zero-mean stochastic processes arises in the 
study of tracking systems \cite{monakov2013physical}. The present paper focuses on studying the multiplicative inversion of non-zero-mean 
complex processes for the specific case where the inverted process has Gaussian statistics. Despite the seemingly simple extension of the 
model, the mathematics involved in deriving the spectral properties of the extended object are challenging. The main result of the paper is 
the derivation of a closed-form expression of the autocorrelation function of such a process, as a function of the first and second order 
statistics of the underlying Gaussian process $1/\mathbf{x}(t)$. The paper is organized as follows: section \ref{sec:definitions} is devoted 
to formalize the problem in terms of studied objects and assumptions; in section \ref{sec:representation} an integral representation of the 
autocorrelation function is derived, allowing all further steps; sections \ref{sec:convergence} and \ref{sec:integrability} tackle the 
important aspect of the existence (in mathematical terms) of autocorrelation function and power spectrum, as well as the possibility to 
establish a power series representation; section \ref{sec:series} is devoted to the derivation of the power series expansion of the 
autocorrelation function; in section \ref{sec:recurrence} a recurrence relation between the coefficients of the previously derived series 
expansion is demonstrated; in section \ref{sec:differential} a first order differential equation is derived from the previously found recurrence 
relation; in section \ref{sec:closedForm} the closed form expression of the autocorrelation function is obtained by integration of the 
differential equation; in section \ref{sec:example} the closed-form formula is applied to a sample case, and finally concluding remarks are summarised in section \ref{sec:conclusions}.

\section{Definitions and Assumptions}
\label{sec:definitions}

Let $\mathbf{w}(t)$ be a zero-mean complex Gaussian process, with real components $\mathbf{a}(t)$, $\mathbf{b}(t)$
\begin{equation}
\label{eq: processW}
\mathbf{w}(t) = \mathbf{a}(t) + i \mathbf{b}(t).
\end{equation}

The process $\mathbf{w}(t)$ is assumed to be Wide Sense Stationary (WSS) \cite{papoulis1965random}, with autocorrelation function defined as
\begin{equation}
\label{eq: autocorrelationW}
R_{ww}(\tau)=E[\mathbf{w}(t+\tau)\mathbf{w}(t)^*],
\end{equation}
where $E[\cdot]$ is the statistical ensemble expectation operator and $(\cdot)^*$ denotes complex conjugation. We define the normalized 
function $r(\tau)$ as follows:
\begin{equation}
\label{eq: defineRandPhi}
r(\tau)\stackrel{\triangle}{=}\frac{R_{ww}(\tau)}{R_{ww}(0)}=|r(\tau)|\exp(i\varphi(\tau)).
\end{equation}
\vspace{1ex}

In view of the above assumptions, it is \cite{papoulis1965random}
\begin{equation}
r(\tau) = \rho(\tau) - i \mu(\tau),
\end{equation}
with
\begin{equation}
\label{eq: defineRhoMu}
\begin{array}{l}
\rho(\tau) \stackrel{\triangle}{=} 2 R_{aa}(\tau) / R_{ww}(0), \\
\mu(\tau) \stackrel{\triangle}{=} 2 R_{ab}(\tau) / R_{ww}(0).
\end{array}
\end{equation}
\vspace{1ex}

The function $\rho(\tau)$ is an even function of $\tau$ with $\rho(0)=1$, whereas $\mu(\tau)$ is an odd function of $\tau$. Furthermore, it is
\begin{equation}
|r(\tau)| \leq r(0) = 1,
\end{equation}
and we also assume that, as for most regular processes mentioned in \cite{papoulis1965random}
\begin{equation}
\lim_{|\tau| \to \infty} r(\tau) = 0.
\end{equation}

The objective of this work is to study the process $\mathbf{s}(t)$ defined as
\begin{equation}
\label{eq: processS}
\mathbf{s}(t) ={ \frac{1}{\mathbf{w}(t) + w_0}},
\end{equation}
where $w_0$ is a complex number.

\section{An integral Representation of the Autocorrelation Function}
\label{sec:representation}

The autocorrelation of $\mathbf{s}(t)$ is formed using Eqs. \eqref{eq: processW}, \eqref{eq: processS} into Eq. \eqref{eq: autocorrelationW}
\begin{equation}
\label{eq: expectionRss}
R_{ss}(\tau) = E \left[ \left( \frac{1}{w_0 + \mathbf{a}(t + \tau) + i\mathbf{b}(t + \tau)} \right) \left( \frac{1}{w_0^* + \mathbf{a}(t) - i\mathbf{b}(t)} \right) \right]
\end{equation}
\vspace{1ex}
where we have anticipated that the process $\mathbf{s}(t)$ is WSS. 

The covariance matrix for $\mathbf{a}(t), \mathbf{b}(t), \mathbf{a}
(t+\tau), \mathbf{b}(t+\tau)$ is, based on definitions \eqref{eq: defineRhoMu}
\begin{equation}
\mathbf{R} = \frac{R_{ww}(0)}{2}
\begin{bmatrix}
1 & 0 & \rho(\tau) & -\mu(\tau) \\
0 & 1 & \mu(\tau) & \rho(\tau) \\
\rho(\tau) & \mu(\tau) & 1 & 0 \\
-\mu(\tau) & \rho(\tau) & 0 & 1
\end{bmatrix}.
\end{equation}

The expectation of Eq. \eqref{eq: expectionRss} is given by
\begin{equation}
\label{eq: RssSimpleIntegral}
R_{ss}(\tau) = \int_{-\infty}^{\infty} d^4 \mathbf{x} \bigg( \frac{p(\mathbf{x})}{(w_0 + x_3 + ix_4)(w_0^* + x_1 - ix_2)} \bigg),
\end{equation}
where $\mathbf{x}=\{x_1, x_2, x_3, x_4\}^T$ and $p(\mathbf{x})$ is a 4-dimensional Gaussian distribution built upon the covariance matrix 
$\mathbf{R}$
\begin{equation}
p(\mathbf{x}) = \frac{1}{(2\pi)^2} \frac{1}{\sqrt{|\mathbf{R}|}} \exp \left( -\frac{1}{2} \mathbf{x}^T \cdot (\mathbf{R}^{-1}) \cdot \mathbf{x} \right).
\end{equation}

In the above definition the dot denotes scalar product, and $\{\cdot\}^T$ indicates vector or matrix transposition. 
By various steps Eq. \eqref{eq: RssSimpleIntegral} becomes
\begin{equation}
\label{eq: RssIntegralRep1}
R_{ss}(\tau) = \frac{1}{R_{ww}(0)} \widehat{R}_{ss}(\tau),
\end{equation}
where we have introduced the normalized autocorrelation function $\widehat{R}_{ss}(\tau)$
\begin{equation}
\label{eq: RssIntegralRep2}
\begin{split}
\widehat{R}_{ss}(\tau) = & -\frac{1}{(2\pi)^2} \int_{0}^{\infty} \exp \left( -\frac{1}{4}(v_1^2 + v_2^2) \right) dv_1 dv_2  \\
& \cdot\int_{0}^{2\pi} \exp \left( -\frac{|r|v_1 v_2}{2} \cos(q_1 - q_2 + \varphi) + i\omega v_1 \cos q_1 + i\omega v_2 \cos q_2 + i(q_1 - q_2) \right) dq_1 d q_2,
\end{split}
\end{equation}
which depends upon the real and non-negative parameter $\omega$ defined as follows:
\begin{equation}
\label{eq: smallomega}
\omega = \frac{|w_0|}{\sqrt{R_{ww}(0)}},
\end{equation}
and where, we recall for convenience, $|r|$ and $\varphi$ were defined in Eq. \eqref{eq: defineRandPhi}. The derivations transforming the 
integral of Eq. \eqref{eq: RssSimpleIntegral} into the Eqs. \eqref{eq: RssIntegralRep1} and \eqref{eq: RssIntegralRep2} include the 
application of a unitary transformation of variables eliminating the dependency from the phase of $w_0$, the representation of the 4-
dimensional Gaussian distribution $p(\cdot)$ through its characteristic function, the application of a linear variable transformation 
eliminating the explicit dependency of the denominator from $|w_0|$, two transformations into double polar coordinates as well as the use of 
Eq. (3.338) in \cite{gradshteyn2007table}, and finally a change of radial variables to normalize the final expression with respect to
$R_{ww}(0)$.

\section{Absolute Convergence and Continuity of the Autocorrelation Function}
\label{sec:convergence}

The above integral \eqref{eq: RssIntegralRep2} is absolutely convergent for $\tau\neq0$ (or equivalently $r\neq 1$), as shown below
\begin{equation}
\label{eq: absoluteConvergence}
\begin{split}
\int_{0}^{\infty} \int_{0}^{2\pi} \bigg| & \exp \left( -\frac{1}{4}(v_1^2 + v_2^2) -\frac{|r|v_1 v_2}{2} \cos(q_1 - q_2 + \varphi) \right) \\
& \cdot \exp \left( i\omega v_1 \cos q_1 + i\omega v_2 \cos q_2 + i(q_1 - q_2) \right) \bigg| dq_1 dq_2 dv_1 dv_2 \\
& < \frac{4\pi^2}{\sqrt{1-|r|^2}} \left( \pi + 2 \tan^{-1} \left( \frac{|r|}{\sqrt{1-|r|^2}} \right) \right).
\end{split}
\end{equation}

The inequality \eqref{eq: absoluteConvergence} can be obtained by distributing the absolute module to the various factors within the integral, 
and by applying Eq. BI (81)(7) and Eq. TI (253), FI II 94 of \cite{gradshteyn2007table}. The right term in the above inequality is finite when 
$|r| < 1$, confirming the absolute convergence of the integral in Eq. \eqref{eq: RssIntegralRep2} for $\tau \neq 0$. Furthermore. the 
integrand in Eq. \eqref{eq: RssIntegralRep2} is visibly a continuous function of the parameter $\omega$, in the domain $\omega \in (0,\infty)$.

\section{Absolute Integrability of the Autocovariance Function}
\label{sec:integrability}

The normalized autocovariance function is obtained as follows
\begin{equation}
\begin{split}
\widehat{C}_{ss}(\tau) = & \widehat{R}_{ss}(\tau) - \lim_{|\tau| \to \infty} \widehat{R}_{ss}(\tau) = -\frac{1}{(2\pi)^2} \int_{0}^{\infty} \exp \left( -\frac{1}{4}(v_1^2 + v_2^2) \right) dv_1 dv_2 \\
& \cdot \int_{0}^{2\pi} \left( \exp \left( -\frac{|r|v_1 v_2}{2} \cos(q_1 - q_2 + \varphi) \right) - 1 \right) \\
& \cdot \exp \left( i\omega (v_1 \cos q_1 + v_2 \cos q_2) + i(q_1 - q_2) \right) dq_1 d q_2.
\end{split}
\end{equation}

We investigate whether $\widehat{C}_{ss}(\tau)$ is absolutely integrable, \textit{i.e.} under which assumptions the following condition 
is fulfilled  
\begin{equation}
\label{eq: conditionAbsoluteIntegrability}
\int_{-\infty}^{\infty} |\widehat{C}_{ss}(\tau)| d\tau < \infty.
\end{equation}

First, we have
\begin{equation}
\label{eq: inequalityCss}
\begin{split}
\int_{-\infty}^{\infty} |\widehat{C}_{ss}(\tau)| d\tau < & \frac{1}{(2\pi)^2} \int_{-\infty}^{\infty} d\tau \int_{0}^{\infty} \exp \left( -\frac{1}{4}(v_1^2+v_2^2) \right) dv_1 dv_2 \\
& \cdot \int_{0}^{2\pi} \bigg| \exp \left( -\frac{|r|v_1 v_2}{2} \cos(q_1-q_2+\varphi) \right) - 1 \bigg| dq_1 dq_2.
\end{split}
\end{equation}

We can find the following expression for the inner integral.
\begin{equation}
\label{eq: Struve}
\int_{0}^{2\pi} \left| \exp \left( -\frac{|r|v_1 v_2}{2} \cos(q_1 - q_2 + \varphi) \right) - 1 \right| dq_1 dq_2 = 4\pi^2 L_0 \left( \frac{|r|v_1 v_2}{2} \right)
\end{equation}
where $L_0(\cdot)$ is the modified Struve function with parameter zero. Eq. \eqref{eq: Struve} is obtained by simple algebraic manipulations 
of the integral, and by using the integral representation of the modified Struve Function 12.2.2 in \cite{abramowitz1965handbook}. When 
introducing Eq. \eqref{eq: Struve} in the inequality \eqref{eq: inequalityCss}, after simple manipulations and by using the series expansion 
of $L_0(\cdot)$ we obtain
\begin{equation}
\label{eq: inequalityCss2}
\int_{-\infty}^{\infty} |\widehat{C}_{ss}(\tau)| d\tau < \int_{-\infty}^{\infty} \frac{4}{\pi} |r| \cdot {}_3F_2 \left( 1,1,1; \frac{3}{2}, \frac{3}{2}; |r|^2 \right) d\tau
\end{equation}
where ${}_3F_2(\cdot)$ is the generalized hypergeometric series defined \textit{e.g.} in section 9.14 of \cite{gradshteyn2007table} with 
$p=3$, $q=2$, zero-balanced according to \cite{evans1984asymptotic}. The above inequality \eqref{eq: inequalityCss2} allows determining 
whether the condition \eqref{eq: conditionAbsoluteIntegrability} is satisfied for specific values of $r$. The improper integral on the right 
side of Eq. \eqref{eq: inequalityCss2} can diverge for two reasons: for the singularity of the integrand at $\tau = 0$, when $r = 1$, and for 
the behaviour of the integrand for $|\tau| \to \infty$, when $|r| \to 0$. As far as the first condition is concerned, the singularity at $\tau 
= 0$ ($|r|=1$) is of logaritmic type \cite{evans1984asymptotic} \cite{rankin1961notebooks}, thus integrable for any typical form of $r$. 
Concerning the behaviour for $|\tau| \to \infty$, the integrand approaches $\frac{4}{\pi} |r|$, and therefore it is sufficient that $|r|$ 
decays faster than $1/|\tau|^\alpha$ with $\alpha >1$ to ensure absolute integrability. As an example, we analyze the case where $r$ has the 
following expression
\begin{equation}
\label{eq: example1}
r = \exp(-a|\tau|), \quad a > 0,
\end{equation}
which corresponds to a Lorentzian power spectrum of the process $\mathbf{w}(t)$. When using Eq. \eqref{eq: example1} into 
Eq. \eqref{eq: inequalityCss2} we obtain
\begin{equation}
\int_{-\infty}^{\infty} |\widehat{C}_{ss}(\tau)| d\tau < \frac{28\zeta(3)}{a\pi} - \frac{8C}{a} \cong \frac{3.39}{a}, \quad a > 0
\end{equation}
where $\zeta(3)$ and $C$ are the Apéry's and Catalan's constants respectively. The above expression is obtained immediately by a change of 
variable in the integral on the right side of Eq. \eqref{eq: inequalityCss2}, and by using the power series definition of the generalized 
hypergeometric function for the final integration. By similar steps, however without obtaining a closed formula, we analyze the following 
normalized autocorrelation of a Gaussian power spectrum
\begin{equation}
r = \exp(-a\tau^2), \quad a > 0,
\end{equation}
and obtain the inequality
\begin{equation}
\int_{-\infty}^{\infty} |\widehat{C}_{ss}(\tau)| d\tau < \frac{4.53}{\sqrt{a}}, \quad a > 0.
\end{equation}

The two examples above show that for certain power spectra of the process $\mathbf{w}(t)$ the normalized autocovariance function 
$\widehat{C}_{ss}(\tau)$ is absolutely integrable, which guarantees, in those cases, the existence of the covariance spectrum of the process 
$\mathbf{s}(t)$. We recall however that the fulfilment of the inequality \eqref{eq: conditionAbsoluteIntegrability} is a sufficient but not 
necessary condition for the existence of the power spectrum, which implies that the power spectrum may exist also when such condition is not 
fulfilled. 

\section{Series Expansion for the Autocorrelation and Autocovariance functions}
\label{sec:series}

It is useful to expand the normalised function $\widehat{R}_{ss}(\tau)$ in power series of $\omega$
\begin{equation}
\label{eq: series1}
\widehat{R}_{ss}(\tau) = \sum_{n=0}^{\infty} \frac{\Omega_n}{n!} \omega^n,
\end{equation}
where $\Omega_n$ is given by
\begin{equation}
\label{eq: omegaAsDerivative}
\Omega_n = \lim_{\omega \to 0} \left( \frac{d^n \widehat{R}_{ss}(\tau)}{d\omega^n} \right).
\end{equation}

By deriving Eq. \eqref{eq: RssIntegralRep2} with respect to $\omega$ and by performing the subsequent limit operations (both operations under 
the integral sign, as justified by the results of the previous sections) the following is immediately obtained
\begin{equation}
\label{eq: omegan}
\begin{split}
\Omega_n = & -\frac{1}{(2\pi)^2} \int_{0}^{\infty} \exp \left( -\frac{1}{4}(v_1^2 + v_2^2) \right) dv_1 dv_2 \cdot \\
& \cdot \int_{0}^{2\pi} \exp \left( -\frac{|r|v_1 v_2}{2} \cos(q_1 - q_2 + \varphi) \right) (i)^n \cdot \\
& \cdot \left[ (v_1 \cos q_1 + v_2 \cos q_2)^n \right] \exp \left( i(q_1 - q_2) \right) dq_1 d q_2.
\end{split}
\end{equation}

It can be shown by simple inequalities and by use of Eq. 3.326.2 in \cite{gradshteyn2007table} that it is
\begin{equation}
|\Omega_n| < 2^{\frac{3n}{2}-1} n\pi \frac{1}{(1 - |r|)^{1+\frac{n}{2}}} \Gamma\left(\frac{n}{2}\right).
\end{equation}

It is easy to verify, \textit{e.g.} by the ratio test, that the series \eqref{eq: series1} converges absolutely for $|r|<1$ as it is
\begin{equation}
\lim_{n \to \infty} \frac{|\Omega_{n+1}| \omega^{n+1} / (n+1)!}{|\Omega_n| \omega^n / n!} = \lim_{n \to \infty} \frac{\Gamma\left(\frac{1+n}{2}\right)}{\Gamma\left(1+\frac{n}{2}\right)} \frac{\sqrt{2}\omega}{\sqrt{1-|r|}} = 0, \quad |r| < 1
\end{equation}

After several manipulations of the Eq. \eqref{eq: omegan} we arrive at the following result
\begin{equation}
\label{eq: omeganFormulaOdd}
[\Omega_n]_{n \text{ odd}} = 0,
\end{equation}
and
\begin{equation}
\label{eq: omeganFormulaEven}
\begin{aligned}
&[\Omega_n]_{n\text{ even}} \\
&= 2 \sum_{k=0}^{n/2-1} \sum_{j=0}^k \frac{(-1)^{(k+\frac{n}{2})} n!}{(n/2-j)!(k-j)!} \left\{ \frac{{}_2\text{F}_1\left(1+j, 1+j-k+\frac{n}{2}, 2+2j-k, |r|^2\right)}{\Gamma(2+2j-k)} \right\} (r^*)^{2j-k+1} \\
&+ \sum_{j=0}^{n/2} \frac{n!}{[(n/2-j)!]^2} \left\{ \frac{{}_2\text{F}_1\left(1+j, 1+j, 2+2j-\frac{n}{2}, |r|^2\right)}{\Gamma\left(2+2j-\frac{n}{2}\right)} \right\} (r^*)^{2j-\frac{n}{2}+1},
\end{aligned}
\end{equation}
where ${}_2F_1(\cdot)$ is the Gauss hypergeometric series (\textit{e.g.} defined in section 15 of \cite{abramowitz1965handbook}) The 
derivation of the above expression from Eq. \eqref{eq: omegan} proceeds by tackling the inner integral over the variables $q_1$, $q_2$ with the 
use of the inegral representation of the modified Bessel function (\textit{e.g.} Eq. WA 201(4) in \cite{gradshteyn2007table}), the identity 
6.631 and the integral formula 8.405 KU 46(1) in \cite{gradshteyn2007table}. After several subsequent manipulations, and by use of the 
integral formula EH I 269(5) in \cite{gradshteyn2007table}, Eqs. \eqref{eq: omeganFormulaOdd}, \eqref{eq: omeganFormulaEven} are obtained. Eq. 
\eqref{eq: omeganFormulaEven} is a finite sum, it allows an easy computation of the coefficient $\Omega_n$ for any order. For example, the 
coefficients $\Omega_n$, computed for $n$ even up to $2$, are reported in Eq. \eqref{eq: omeganExampleComplex} for the general case where $r$ 
is complex:
\begin{equation}
\label{eq: omeganExampleComplex}
\begin{array}{rll}
\Omega_0 & = \displaystyle -\frac{\ln(1 - |r|^2)}{r} & \\
\Omega_2 & = \displaystyle \frac{2(1 - 2r^* + r^*/r)}{1 - |r|^2} + \frac{2\ln(1 - |r|^2)}{r^2}, & r \in \mathbb{C} \\
\Omega_4 & = \dots & 
\end{array}
\end{equation}

In the case that $r$ is real (symmetric power spectrum of $\mathbf{w}(t)$) the expression becomes simpler, below the coefficients $\Omega_n$ 
with $n$ up to $6$ are reported:
\begin{equation}
\label{eq: omeganExampleReal}
\begin{array}{rll}
\Omega_0 & = \displaystyle -\frac{\ln(1 - r^2)}{r} \\[3ex]
\Omega_2 & = \displaystyle \frac{4}{1 + r} + \frac{2\ln(1 - r^2)}{r^2} \\[3ex]
\Omega_4 & = \displaystyle -\frac{12}{r} - \frac{24}{(1 + r)^2} - \frac{12\ln(1 - r^2)}{r^3} & r \in \mathbb{R} \\[3ex]
\Omega_6 & = \displaystyle \frac{120}{r^2} + \frac{320}{(1 + r)^3} + \frac{80}{(1 + r)^2} + \frac{80}{1 + r} + \frac{120\ln(1 - r^2)}{r^4} \\[3ex]
\Omega_8 & = \dots
\end{array}
\end{equation}

We can also compute the limit expression of the coefficients $\Omega_n$ when $\tau \to \infty$ (and therefore when $|r| \to 0$). The only 
terms that survive in Eq. \eqref{eq: omeganFormulaEven} in such a limit are those for which $k$ is odd and $2j-k+1=0$, which immediately leads 
to 
\begin{equation}
\lim_{\tau \to \infty} [\Omega_n]_{n \text{ even}} = 2(-1)^{n/2 + 1} \frac{(2^{n/2} - 1)n!}{(n/2 + 1)!}
\end{equation}
and in turn to
\begin{equation}
\label{eq: limitRss}
\lim_{\tau \to \infty} \widehat{R}_{ss}(\tau) =  \frac{(1 - \exp(-\omega^2))^2}{\omega^2}.
\end{equation}

The above coincides, as expected, with the absolute square of the mean of the process $\mathbf{s}(t)$ (normalized by $R_{ww}(0)$, which can be 
directly and easily computed as $|E[\mathbf{s}(t)]|^2$ (not done here). Based on the above discussion, the covariance of the process 
$\mathbf{s}(t)$ can be expanded as follows:
\begin{equation}
\widehat{C}_{ss}(\tau) = \sum_{n=0}^{\infty} \frac{\Omega'_n}{n!} \omega^n,
\end{equation}
where
\begin{equation}
\label{eq: omeganPrime}
\Omega'_n = \Omega_n - \lim_{\tau \to \infty} \Omega_n.
\end{equation}

\section{Recurrence relation}
\label{sec:recurrence}

In this section we prove that the following recurrence relation holds for the coefficients $\Omega_n$ defined in the previous section
\begin{equation}
\label{eq: recursion}
\Omega_{2(m+1)} r + 2(1 + 2m)\Omega_{2m} = \frac{(-1)^m 2^{3+2m} \Gamma\left(\frac{3}{2} + m\right)}{(1 + m)\sqrt{\pi}} (1 - 2^m a^{m+1}), m \ge 0
\end{equation}
where $a$ is defined as follows
\begin{equation}
\label{eq: aDefinition}
a = \frac{1 - r \text{e}(r)}{1 - |r|^2}.
\end{equation}

Although the above recurrence relation was established “experimentally” by use of Eq. \eqref{eq: omeganFormulaEven}, by inspecting expressions 
of $\Omega_n$ evaluated for various (up to large) even orders $n$, a formal proof based on the use of Eq. \eqref{eq: omeganFormulaEven} could 
not be obtained. Instead it was possible to prove the recurrence relation by the steps illustrated below. Firstly, the integral of Eq. 
\eqref{eq: RssIntegralRep2} is transformed into the following integral
\begin{equation}
\label{eq: integralRepresenationTrigs}
\widehat{R}_{ss}(\tau) = -\frac{1}{2\pi} \int_{0}^{\pi} d\theta \int_{0}^{2\pi} e^{i\phi} \frac{\exp\left[ -\omega^2 \frac{1 + \sin\theta \cos\phi}{1 + \sin\theta (\text{re}(r) \cos\phi - \text{im}(r) \sin\phi)} \right]}{1 + \sin\theta (\text{re}(r) \cos\phi - \text{im}(r) \sin\phi)} d\phi.
\end{equation}

The above expression is obtained by first normalizing through $\omega$ the radial variables in Eq. \eqref{eq: RssIntegralRep2} followed by 
integration with respect to the angular variables $q_1$, $q_2$ by use of the integral formula 3.915 2. in \cite{gradshteyn2007table}. After 
changing the residual radial variables to polar variables, the resulting integral is performed with the help of Eq. 6.631 1. from 
\cite{gradshteyn2007table}, followed by simple rescaling of variables. When expanding the exponential term in Eq. \eqref{eq: 
integralRepresenationTrigs} in power series of $\omega^2$, we obtain the following expression of the coefficient $\Omega_n$ already defined in 
Eq. \eqref{eq: omegaAsDerivative}
\begin{equation}
\label{eq: omega2m}
\Omega_{2m} = -\frac{1}{\pi} \frac{(-1)^m (2m)!}{m!} \int_{0}^{\pi/2} d\theta \int_{0}^{2\pi} e^{i\phi} \frac{A^m}{B^{m+1}} d\phi
\end{equation}
having defined
\begin{equation}
\label{eq: defAB}
\begin{aligned}
A &= 1 + \sin\theta \cos\phi \\
B &= 1 + \sin\theta (\text{re}(r) \cos\phi - \text{im}(r) \sin\phi)
\end{aligned}
\end{equation}

When introducing the Eq. \eqref{eq: omega2m} into Eq. \eqref{eq: recursion}, and after various algebraic manipulations, we can transform the 
to-be-proven relation \eqref{eq: recursion} into the following complex equality to be proven for every integer $m\geq0$
\begin{equation}
\label{eq: recursion2}
\int_{0}^{\pi/2} d\theta \int_{0}^{2\pi} e^{i\phi} \left( r \frac{A^{m+1}}{B^{m+2}} - \frac{A^m}{B^{m+1}} \right) d\phi = \frac{2\pi}{(1 + m)} (1 - 2^m a^{m+1})
\end{equation}

As the right end side of the above equation is purely real, the imaginary component of the integral must be zero. This fact can be easily 
proved by observing that
\begin{equation}
\text{im}\left[ e^{i\phi} \left( r \frac{A^{m+1}}{B^{m+2}} - \frac{A^m}{B^{m+1}} \right) \right] = \left(\frac{A}{B}\right)^m \frac{1}{\sin\theta} \frac{\partial}{\partial\phi} \left(\frac{A}{B}\right)
\end{equation}
which, in view of the expression of $A$, $B$ in Eq. \eqref{eq: defAB} shows that the inner $\phi$-integral in Eq. \eqref{eq: recursion2} is 
identically zero. We are then left with proving the equality for real components of Eq. \eqref{eq: recursion2} which is easily converted into 
the following:
\begin{equation}
\label{eq: recursion3}
\int_{0}^{\pi/2} d\theta \int_{0}^{2\pi} \left(\frac{A}{B}\right)^m \frac{C}{B^2} d\phi = \frac{2\pi}{(1 + m)} (1 - 2^m a^{m+1})
\end{equation}
with
\begin{equation}
C = (-1 + \text{re}(r)) \cos\phi - \text{im}(r) \sin\phi
\end{equation}

The above Eq. \eqref{eq: recursion3} can be proven through the method of generating function. Each side of Eq. \eqref{eq: recursion3} is 
multiplied by $(1+m)t^m$, followed by summing over $m$ from $0$ to $\infty$. After few manipulations the following equation is obtained
\begin{equation}
\label{eq: recursion4}
\int_{0}^{\pi/2} d\theta \int_{0}^{2\pi} \frac{C}{(B - At)^2} d\phi = \frac{2\pi}{(1 - t)} - \frac{2\pi a}{(1 - 2at)}
\end{equation}

The above equation can be verified by using Eq. 3.661 4. of \cite{gradshteyn2007table} (properly manipulated to serve our purposes, including 
derivation with respect to its parameters) to perform the inner integral, leading to
\begin{equation}
\label{eq: recursion5}
\int_{0}^{\pi/2} \left( -\frac{2\pi \sin\theta \left( (\text{re}(r) - t)(-1 + \text{re}(r)) + \text{im}(r)^2 \right)}{\left( (1 - t)^2 - \sin^2\theta \left( (\text{re}(r) - t)^2 + \text{im}(r)^2 \right) \right)^{3/2}} \right) d\theta = \frac{2\pi}{(1 - t)} - \frac{2\pi a}{(1 - 2at)}
\end{equation}

The proof of the above equality is obtained by the variable substitution $z=\sin\theta$ followed by the use of Eq. 2.271 5. of 
\cite{gradshteyn2007table}. The coincidence of the two sides in Eq. \eqref{eq: recursion5} implies that Eq. \eqref{eq: recursion3} is true for 
all $m\geq0$, which in turn ensures the same for the initial recurrence relation \eqref{eq: recursion}.
\section{First order differential equation}
\label{sec:differential}

The recurrence relation \eqref{eq: recursion} proven in the previous section leads to the following differential equation for the normalized 
autocorrelation function $\hat{R}_{ss}$
\begin{equation}
\label{eq: diffEquation}
\omega r \frac{d\hat{R}_{ss}}{d\omega} + 2\omega^2 \hat{R}_{ss} = 2\left(1 - 2e^{-\omega^2} + e^{-2a\omega^2}\right)
\end{equation}
where $r$, $\omega$ and $a$ were defined in Eqs. \eqref{eq: defineRandPhi}, \eqref{eq: smallomega} and \eqref{eq: aDefinition}. The above 
equation is obtained by first multiplying each side of Eq. \eqref{eq: recursion} by $\omega^{2m} / (1 + 2m)!$ leading to
\begin{equation}
\frac{1}{\omega} \frac{d}{d\omega} \Omega_{2(m+1)} r \frac{\omega^{2m+2}}{(2m+2)!} + 2\Omega_{2m} \frac{\omega^{2m}}{(2m)!} = \frac{(-1)^m 2^{3+2m} \Gamma\left(\frac{3}{2} + m\right)}{(1+m)\sqrt{\pi}} (1 - 2^m a^{m+1}) \frac{\omega^{2m}}{(2m+1)!}
\end{equation}

Each side of the above equation is summed over $m$ from 0 to $\infty$. After re-indexing of the sum on the left side, and simple algebraic 
manipulations, the differential equation \eqref{eq: diffEquation} is obtained. As for any first order differential equation, a boundary 
condition is required in order to find a unique solution, which in our case is given by Eqs. \eqref{eq: series1} and \eqref{eq: omeganExampleComplex}
\begin{equation}
\label{eq: boundaryCondition}
\hat{R}_{ss}(\tau)\Big|_{\omega=0} = \Omega_0 = -\frac{\ln(1 - |r|^2)}{r}
\end{equation}

\section{Closed form expression of the autocorrelation function}
\label{sec:closedForm}

In the present section we use $\gamma=\omega^2$ to express the differential equation and the consequent solution for the normalised 
autocorrelation function $\hat{R}_{ss}$. In terms of $\gamma$ the differential equation \eqref{eq: diffEquation} is immediately re-written as
\begin{equation} 
\label{eq: diffEquationGamma}
r \frac{d\hat{R}_{ss}}{d\gamma} + \hat{R}_{ss} = \frac{1 - 2e^{-\gamma} + e^{-2a\gamma}}{\gamma}
\end{equation}

The above equation can be easily integrated leading to the following 
\begin{equation}
\label{eq: closedForm}
\widehat{R}_{ss} = 
\begin{cases}
\displaystyle \frac{1}{r} e^{-\frac{\gamma}{r}} \left( \operatorname{Ei}\left(\gamma \frac{1}{r}\right) - 2 \operatorname{Ei}\left(\gamma \frac{1-r}{r}\right) + \operatorname{Ei}\left(\gamma \frac{(1-r)^2}{r(1-|r|^2)}\right) \right) & r \neq 0 \\[12pt]
\displaystyle \frac{(1 - e^{-\gamma})^2}{\gamma} & r = 0
\end{cases}
,
\end{equation}
where $\operatorname{Ei}(\cdot)$ is the standard complex analytic continuation of the Exponential Integral function defined in 5.1.2 of 
\cite{abramowitz1965handbook}. Even though $\operatorname{Ei}(\cdot)$ possesses a branch cut along $(-\infty, 0]$, the linear combination in 
the top expression of Eq. \eqref{eq: closedForm} cancels out the individual discontinuities, leaving the total expression completely free of 
a branch cut. Furthermore, the function has a removable singularity at $r=0$ where it approaches the limit reported in the bottom of Eq. \eqref{eq: closedForm}. In the limit where $\gamma \to 0$, the expression converges to the boundary condition \eqref{eq: boundaryCondition}.

\section{Sample computation}
\label{sec:example}
As an example of use of the result from the previous section, we compute the autocorrelation $R_{ss}$ of the process $\mathbf{s}(t)$ of Eq. \eqref{eq: processS}  in the case that the underlying complex Gaussian process $\mathbf{w}(t)$ has a flat spectrum within a given bandwidth.  
Willing to exercise the case of complex $r$ in Eq. \eqref{eq: closedForm}, we assume that the flat spectrum is offset with respect to  zero 
frequency, so that an imaginary component exists for $r$. The results of the calculation are shown in Fig. 1, together with outputs from 
numerical simulations. The agreement is excellent.

\begin{figure}[htbp]
    \centering
    \begin{subfigure}[b]{0.48\textwidth}
        \centering
        \includegraphics[width=\textwidth]{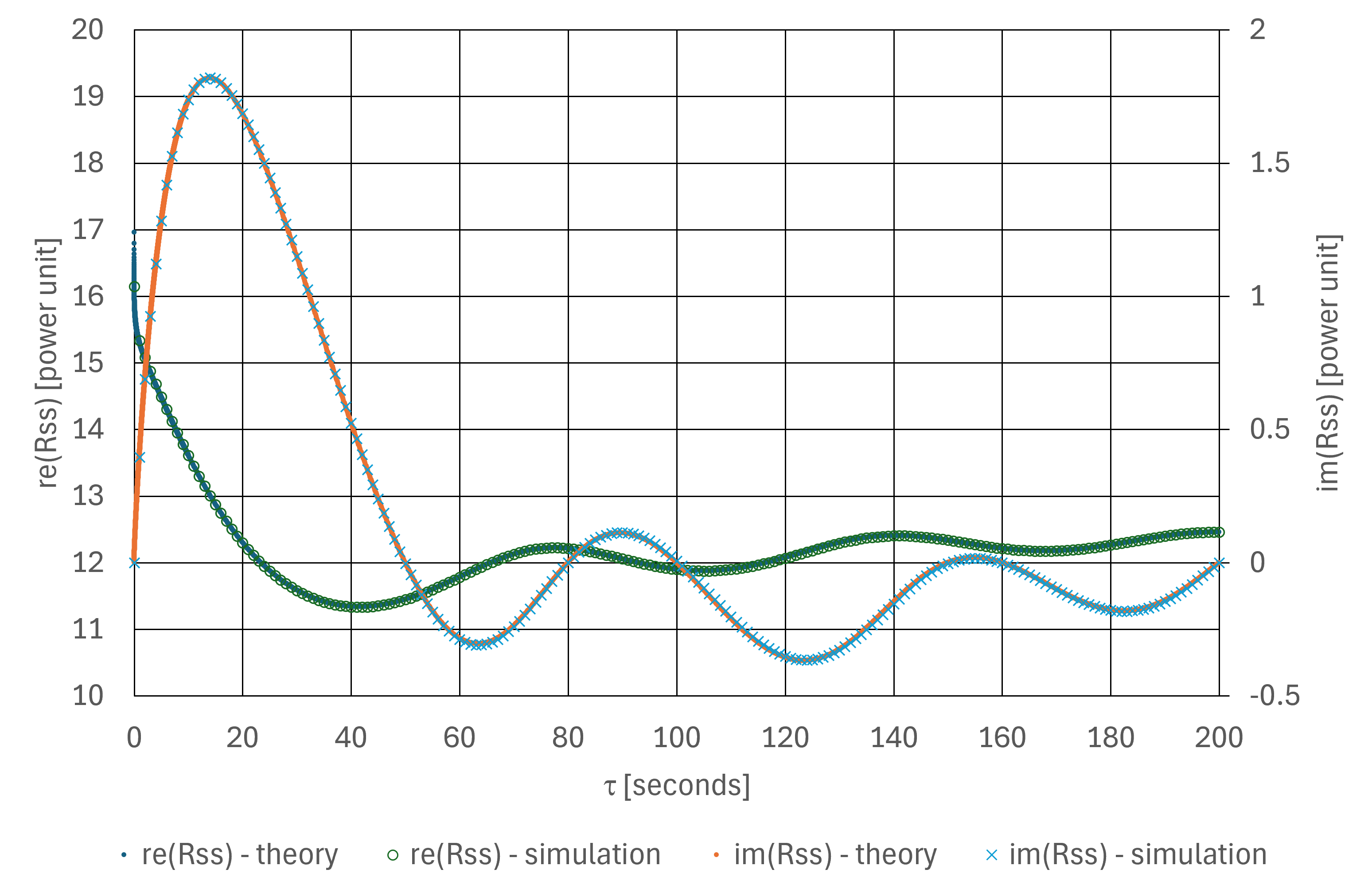}
        \label{fig:1a}
    \end{subfigure}
    \hfill
    \begin{subfigure}[b]{0.48\textwidth}
        \centering
        \includegraphics[width=\textwidth]{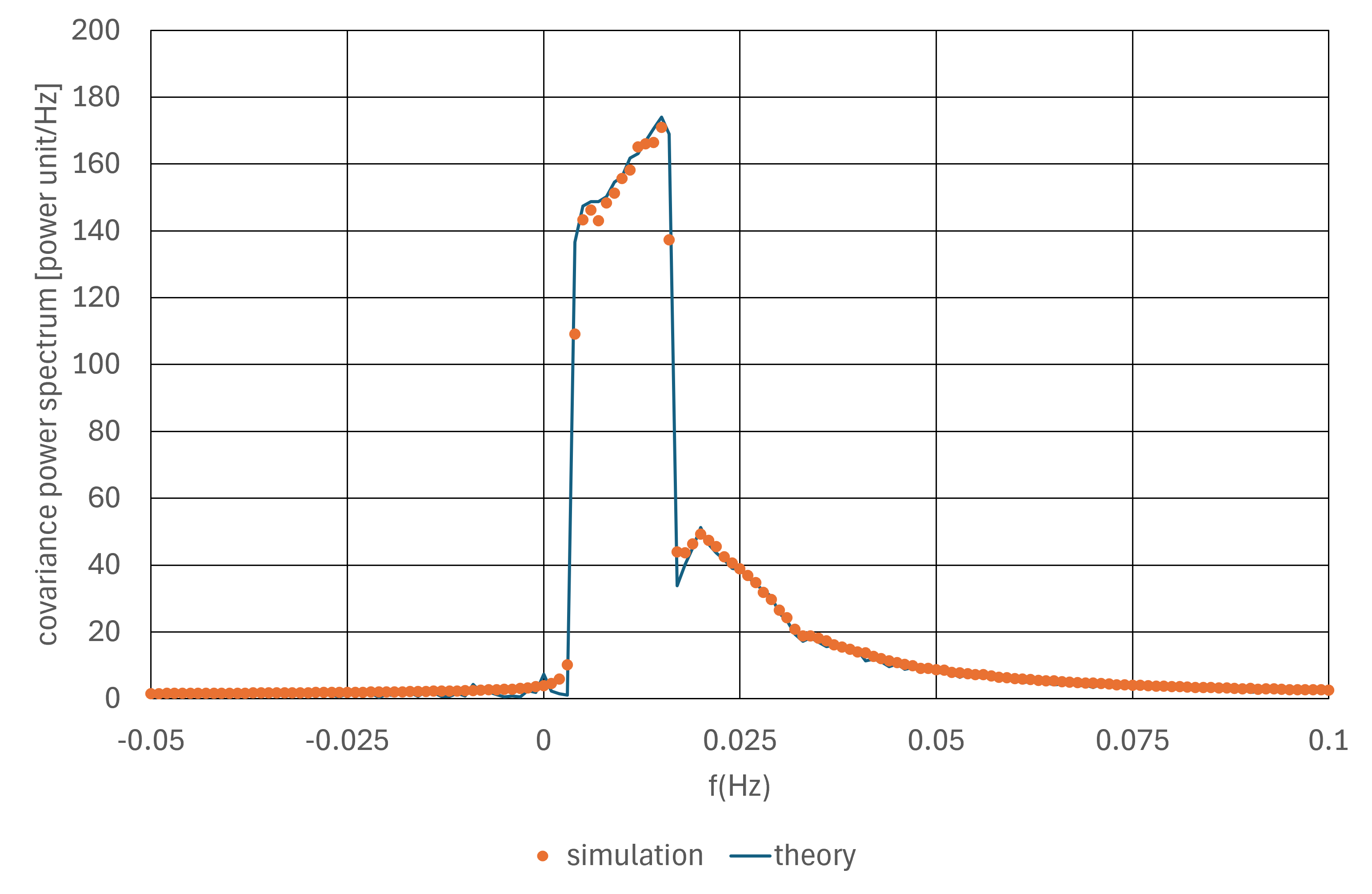}
        \label{fig:1b}
    \end{subfigure}
    \caption{sample computation of the autocorrelation $R_{ss}$ of the process $\mathbf{s}(t)$ defined in Eq. \eqref{eq: processS}  through
    Eqs. \eqref{eq: closedForm} and \eqref{eq: RssIntegralRep1}. The considered power spectrum of the underlying complex Gaussian process
    $\mathbf{w}(t)$ has unitary flat spectral density centred at $1/100$ Hz, with bandwidth of $1/80$ Hz. The value of $w_0$ was set to the arbitrary
    value of $2/7 \exp(\mathrm{i} \pi / 7)$. On the left the real and imaginary parts of $R_{ss}$ are shown, compared with results from 
    numerical simulations. On the right the resulting covariance power spectrum is shown, also compared with results from simulations.
    Care must be taken in computing $R_{ss}$ through Eq. \eqref{eq: closedForm} for $\tau \to 0$ where $r \to 1$ due to the 
    divergence to infinity (as hinted in the plot of the real part) and for $\tau \to \infty$ ($r \to 0$) due to the
    removable singularity. Furthermore, due to the unbound spectral extension of the process $\mathbf{s}(t)$, care must be taken
    to ensure that the simulation approaches a time-continuous process (sufficiently narrow bandwidth vs. sampling rate).}
    \label{fig:1}
\end{figure}

\section{Conclusions}
\label{sec:conclusions}
We have computed the autocorrelation of the multiplicative inverse of a non zero mean complex Gaussian process. The key results are the closed-form expression of the autocorrelation given in Eq. \eqref{eq: closedForm} and the first order differential equation to which such function obeys, given in Eq. \eqref{eq: diffEquationGamma}.

\bibliographystyle{unsrt}
\bibliography{references}

@book{deutsch2017nonlinear,
  title={Nonlinear transformations of random processes},
  author={Deutsch, Ralph},
  year={2017},
  publisher={Courier Dover Publications}
}

@article{monakov2013physical,
  title={Physical and statistical properties of the complex monopulse ratio},
  author={Monakov, Andrei},
  journal={IEEE Transactions on Aerospace and Electronic Systems},
  volume={49},
  number={2},
  pages={960--968},
  year={2013},
  publisher={IEEE}
}

@book{papoulis1965random,
  title={Random variables and stochastic processes},
  author={Papoulis, Athanasios},
  year={1965},
  publisher={McGraw Hill}
}

@article{gradshteyn2007table,
  title={Table of integrals, series, and products. Seventh},
  author={Gradshteyn, Israil Solomonowitsch and Ryzhik, IM},
  journal={Elsevier/Academic Press, Amsterdam},
  volume={48},
  pages={1171},
  year={2007}
}

@book{abramowitz1965handbook,
  title={Handbook of mathematical functions: with formulas, graphs, and mathematical tables},
  author={Abramowitz, Milton and Stegun, Irene A},
  volume={55},
  year={1965},
  publisher={Courier Corporation}
}

@article{evans1984asymptotic,
  title={Asymptotic formulas for zero-balanced hypergeometric series},
  author={Evans, Ronald J and Stanton, Dennis},
  journal={SIAM journal on mathematical analysis},
  volume={15},
  number={5},
  pages={1010--1020},
  year={1984},
  publisher={SIAM}
}

@article{rankin1961notebooks,
  title={Notebooks of Srinivasa Ramanujan, Tata Institute of Fundamental Research, Bombay. 2 vols. 4 to boxed. 1957 100 Rs. net.},
  author={Rankin, RA},
  journal={The Mathematical Gazette},
  volume={45},
  number={351},
  pages={73--74},
  year={1961},
  publisher={Cambridge University Press}
}

\end{document}